\documentclass[%
 reprint,
 amsmath,amssymb,
 aps,
]{revtex4-2}
\DeclareUnicodeCharacter{2061}{}
\usepackage[labelformat=simple]{subcaption}
\usepackage{graphicx}
\usepackage{dcolumn}
\usepackage{bm}
\usepackage{amsmath,amssymb}
\usepackage{appendix}
\usepackage{caption}
\usepackage[dvipsnames]{xcolor}

\usepackage[utf8]{inputenc}
\usepackage{mathtools}
\usepackage{url}

\begin{document}

\preprint{APS/123-QED}

\title{Stochastic Non-Linear Influence in Synchronisation Dynamics}

\author{Hannah Gallant}
\affiliation{%
 Defence Science and Technology Group, Canberra 2600, Australia
}%
\affiliation{
 Australian National University, Canberra 2601, Australia
}%


\author{Dale Roberts}
\affiliation{
 Australian National University, Canberra 2601, Australia
}%

\author{Alexander Kalloniatis}
\affiliation{%
 Defence Science and Technology Group, Canberra 2600, Australia
}%


\date{\today}

\begin{abstract}
We propose a mathematical model that represents the influence of an external actor (influencer) on a network of actors (influenced). The model is an adaptation of control systems within the family of formulations inspired by the Kuramoto Model of synchronisation. Capturing influence as a capacity to affect the character or behaviour of another, we study an external node’s ability to influence the synchronisation of the Kuramoto system while its global behaviour is pulled to the external node’s frequency and away from its natural mean frequency. In our work, stochastically generated dynamical weights assigned to the links between the network and the external influencer whereby links from one system to the other are assigned via one-sided heavy tail noise, generated by the Gamma distribution.  We perform numerical experiments to examine transition points in the ability of the external node to alter the behaviour of the influenced network; either to disrupt synchronisation or to drive it to collective frequencies determined by the external node. We examine the dependence of transition points on the external node's frequency, where too ambitious a driving frequency fails to influence the system while retaining a synchronised state, and reducing achieves a state of synchronisation at a frequency shifted from the mean natural frequency. We also look at the analytic approximation for the system close to synchronisation and fragmentation to understand its behaviour around this limit. 
\end{abstract}
\maketitle


\section{\label{sec:intro}Introduction}
The modern phenomenon of social media and its impact on societal behaviours suggests there is value in research on the mechanisms of `influence' within networked systems of interacting agents. This paper proposes a mathematical model of such a system, where synchronisation is the underlying dynamic occurring on a network. 

Influence can be defined as the capacity of an agent or set of agents to have an effect on the thoughts, character, development or behaviour of an individual; clearly such agents (or `entities') require a level of cognitive ability. The concept of influence commonly arises in studies of social networks and collective behaviour \cite{Barrat2008}. West et al discuss social influence within complex social systems as being underpinned by peoples' need to understand and predict their environment, a sense of social identity and belonging, as well as autonomy in regards to decisions and interactions \cite{West2025,West2023}. We draw on aspects elucidated in this work to build out a model of social influence. One such aspect is the dimension of authority-based influence where an authority (either social or formal) shares some sense of identity with a group; \textit{this form of influence is more effective when the authority is seen as representative of the group and their needs, when the influencer is seen within the group being influenced}. Furthermore, West et al make a distinction between autonomy and control, whereby the influencer providing for autonomous motivations is more likely to result in the desired actions from the influenced: \textit{influence effects that are not `heavy-handed' and respect the internal interactions within a group are likely to be more effective} \cite{West2025}. These two features are pivotal for the model of influence we propose and test as it reflects behaviours in complex social groups recognised by psychologists and sociologists.

We set our model of influence within the framework of representations of synchronisation dynamics. (See \cite{Flache2017} for review of classes of models for social influence) The synchronisation of entities in a complex system is pivotal in many scientific areas of application including biological, chemical, physical, and -- for this paper in particular -- social systems \cite{Manrubia2004}. Examples of applications to social and human-technical systems include human-robot interactions \cite{Mizumoto2010}, the synchronisation of pedestrians \cite{Strogatz2005}, rhythmic applause amongst a group \cite{Neda2000}, coherence in opinion dynamics \cite{Pluchino2006} and distributed decision-making \cite{Kalloniatis2020}. Specifically the collective behaviours, like synchronisation, can capture aspects of a social system such as interpersonal synchronisation, opinion dynamics, consensus formation and decision coherence. 
A general description of synchronisation is via a set of coupled oscillators, described by a set of differential equations
\begin{equation}
\dot{\varphi}_{i}=f_{i}(\{\varphi\})
\label{eq:gensynch}
\end{equation}
where $\varphi_{i}$ is the $i^{th}$ oscillator of the system, $f_{i}$ is some function, and $\{\varphi_{i}\}$ is the set of all oscillators. Here, the network can be viewed through the interaction (or dependence) of each oscillator on others within the system. The solution to this general description has the potential to describe a range of types of synchronisation, including phase synchronisation and frequency synchronisation, depending on the form of $f_i$. 

For this paper, we consider the Kuramoto model \cite{Kuramoto1984, Kuramoto2026}, a particularly powerful model of synchronisation, which describes the synchronisation of a large set of $N$ coupled oscillators, $\varphi_{i}$, each with their own intrinsic frequency, $\omega_{i}$. In this model, Eq.~(\ref{eq:gensynch}) takes the form 
\begin{equation}
\dot{\varphi}_{i}=\omega_{i}-\sigma\sum_{j=1}^{N} A_{ij} \sin{\left(⁡ \varphi_{i} - \varphi_{j} \right)} 		
\label{eq:kuramoto}
\end{equation}
where $\sigma$ is the coupling strength between oscillators, the adjacency matrix $A_{ij}$ represents the connections on the network. (see \cite{Strogatz2000,Acebron2005,Rodrigues2016} for reviews). Specifically, each oscillator here represents the decision cycle for a single agent, representing the repeating dynamic nature of decision making through perception of the external environment, decisions of actions to take and then the actions themselves. The important temporal aspects of decision-making including the dynamic nature of thought processes, the impact they have on choices and preferences, as well as the role of expertise are described in \cite{Ariely2001,Herbig2009}. In considering the collective synchronisation here, we thus model the interdependency of this process within a social system. Hence  the synchronisation of the Kuramoto system represents the coherence of the social system, capturing the notion that decision making is an evolving process that is informed by a social network. In this paper we construct a model that represents the external influence of the system's collective {\it synchronisation} behaviour, but where we add a mechanism reflective of our extraction of properties of social influence from \cite{West2025}. 

We achieve a representation of influence in a social system involving synchronisation dynamics by subjecting the oscillator system to an external influence node connected in a random, time-varying way. We can consider this via a random rewiring of the relevant links between the external node and the system. We refer to these as influence links. While synchronisation in systems with random and asymmetric couplings tends to chaos \cite{Sompolinksy1988}, synchronisation in time-varying systems has been widely studied (see \cite{Ghosh2022} for review on synchronisation in time-varying networks). Work in \cite{So2008} demonstrated synchronisation in a network whose connectivity switches value with a given frequency, while identifying mechanisms for both synchronisation and desynchronisation. Recent work in \cite{Groisman2023}, proved the conditions for global synchronisation for a Kuramoto system on a randomly varying graph and \cite{Faggian2019} demonstrated that for such a time-varying random network, synchronisation can still occur even when there is vanishing connectivity in the system. Synchronisation has also been seen to enhance and be more efficient in time-varying systems \cite{Leander2015}, as well as when the network changes under a dynamical, adaptive process  \cite{Zhou2006,Eom2016}. Although this has been considered for a single Kuramoto system, we are particularly interested in adapting similar mechanisms for influence drawing from applications of control in Kuramoto systems.
 
The control problem, such as that solved by the Kalman filter \cite{Stengel1994}, is whether a system can be controlled such that any final state is reachable by any initial state. A network version of this has been developed in \cite{Liu2011} (See \cite{D'Souza2023,Liu2016} for reviews of control in complex systems). For network synchronisation this then becomes whether synchrony can be achieved to a desired frequency $\Omega$, from any initial oscillator state. While some control methods are shown to desynchronise the system \cite{Gjata2017}, previous work by \cite{Brede2016,Kalloniatis2019} has shown that for deterministic dynamics, adaptive mechanisms can be more efficient than traditional pinning control approaches, such as \cite{Skardal2015}. (see \cite{Sorrentino2007,porfiri2008} for references on pinning control in complex networks) 
Here we extend this concept of controlling of a system to influence, whereby the influencing agent seeks to shape or direct the system away from its natural behaviour. In contrast to `control', which is often seen as a pervasive persistent interaction between one agent and another, influence is more subtle, nuanced and intermittent based on the respect for autonomy emphasised by \cite{West2025}. We shall achieve this by a stochastic interaction, to reproduce the non-persistent aspects. Castellano et al states ``Randomness is a necessary ingredient of social interactions: both our individual attitudes and the social influence of our peers may vary in a non-predictable way. Besides, the influence of external factors such as mass media, propaganda, etc., is also hardly predictable. In this respect, opinion dynamics is a stochastic process." \cite{Castellano2009} We represent this inherent time-dependence and stochasticity associated with social systems for our model via the influence links. There is a large body of work modelling non-Gaussian stochastic dynamics on networks (see for example \cite{Tanaka2020,Wang2021}). Specifically for us, we capture the autonomy-respecting nature of influence by using skewed heavy tail noise to reflect the characteristic of random but gentle nudging as a manifestation of influence by one agent over another. For a system that would otherwise synchronise, an example is where the influencing effect seeks through such nudging to either speed up or slow down the collective behaviour of the influenced system away from its natural preferred collective frequency, namely $\bar{\omega}$. In the following, we examine the thresholds by which the influencing node can {\it capture} the system to a preferred frequency $\Omega\neq \bar\omega$, and where collective synchronisation may be diminished {\it and thus the group identity undesirably disrupted} \cite{West2025}. 

\section{\label{sec:model}The Model: Non-Linear Stochastic Influence}
We present a model that represents a network of oscillators, connected via random links to an external oscillator and study how this external node influences the global behaviour of the system. First, consider a network $A_{ij}$ of $N$ oscillators, $\varphi_{i}$, with an external node $\theta$ connected to $\varphi_{i}$ via $B_{i}$. The dynamics of this system are governed by the equations
\begin{align} 
\dot{\theta} = \Omega& \nonumber \\ 
\dot{\varphi}_{i} = \omega_{i}& - \frac{\sigma}{k_{i}} \sum^{N}_{j=1} A_{ij} \sin{ \left(\varphi_{i}-\varphi_{j}\right)} \nonumber\\
& - \tau B_{i} \sin{ \left(\varphi_{i}-\theta\right)}
\label{eq:model}
\end{align}
where $\omega_{i}$ is the natural frequency of node $\varphi_{i}$, $k_i$ is the degree of node $i$, $\Omega$ is the natural frequency of the external node, $\sigma$ is the internal coupling coefficient of the network $A_{ij}$ and $\tau$ is the strength of the links between ${\varphi_{i}}$ and $\theta$. The internal coupling is degree-normalised: each oscillator averages over its neighbours, so $\sigma$ directly controls the coupling strength independently of $N$. We also note that an oscillator $\varphi_{i}$ is connected to the external node when $B_{i}\neq0$, and is not connected for $B_{i}=0$ when it will synchronise according to its adjacent nodes defined by $A_{ij}$. We note here that while $A_{ij}$ represents an undirected link in the internal system, $B_i$ is a directed link in order to represent the asymmetric nature of influence captured in this model. Since $B_i$ does not include a $1/k_i$ factor, $\tau$ has a consistent physical meaning across different topologies: $\tau \sim \sigma$ gives an influence comparable to the internal coupling.  

This formulation can be seen to be consistent with the multiplicative stochastic approach of \cite{KallZup2013}. However, because we treat the noise like a weight we seek a formulation for $B_i$ that is bounded on one side, in contrast to the L{\'e}vy noise used in \cite{KallRob2017, RobKall2018}. The multiplicative noise is introduced in the coupling between the external influencing agent and the system being influenced.  Specifically, to capture the characteristic of subtle nudging as a modifier of a coupling weight on a network link,  including zero values, we use a Gamma distribution bounded from below by zero. The random variable $B_{i}$ is constructed such that this random link is adaptive, and we define the stochastic coupling strength as a random variable via a Gamma distribution \begin{equation}
B_{i}=\frac{G_{i}}{(\varphi_i-\theta)^{2}+1}
\end{equation}
where $G_{i}\sim\Gamma \left(\alpha, \frac{1}{\alpha} \right)$,
with shape parameter $\alpha$ and 
$\mathbb{E}[G_{i}]=1$, is a Lorentzian kernel independent of both $\alpha$ and the graph topology. Smaller choices of $\alpha$ give a peak at or closer to zero; draws of link strengths from such distributions will represent frequent soft nudges of the system. Thus in the model where $\alpha\leq1$, the influencer node guides the otherwise synchronising system through more frequent random weak nudges tending to infrequent strong nudges. This distribution also leads to the expectation value $\mathbb{E}[B_{i}]=\frac{1}{(\varphi_i-\theta)^{2}+1}$, noting that $\mathbb{E}[B_{i}]\leq 1$ where $\mathbb{E}[B_{i}]\approx 1$ when $\varphi_{i} \approx \theta$. We choose this model as we seek to implement a stronger effect between $\theta$ and $\varphi_{i}$ when they are close in phase. Notes that when they are close in phase the form of the interaction will pull phases closer together. If phase slipping occurs and phases jump $2n \pi$ away from each other the interaction is clearly suppressed. The parameter $\alpha$ controls only the variability of the influence link ($\mathbb{CV}[G_i] = \frac{1}{\sqrt{\alpha}}$), while $\tau$ directly sets the mean influence strength independently of the network structure. We note here that we randomly select $B_{i}$ from this distribution at a given time increment in order to create an effect that is reminiscent of blinking but where the influence links while predominantly extremely weak are frequently non-vanishing. $B_i$ can be consider to have its own stochastic time-profile and in doing so, this randomly changing link replicates the stochastic element of these interactions with the social system. 

In considering influence and group identity, we look to those social networks formed by meaningful and reciprocal relationships. The cognitive capacity of people means that such active social networks have a mean size of 150, namely Dunbar's number \cite{Dunbar2018}. This scale for meaningful relationships is also seen to translate to online communities, such as those in social media \cite{Dunbar2015} and online gaming communities \cite{Fuchs2014}, despite the internet allowing an individual to form connections in systems of a much larger size. Human societies, including online communities, are structured social systems where alliances are formed through friendship and emotional closeness. 
Thus these social structures are ultimately founded on a common set of values and convictions \cite{Dunbar2018}. When it comes to social complexities, such as opinion dynamics, clustering and social cohesion, the cliques represented by groups smaller than Dunbar's number would have tightly shared bond.  While ``for groups of sizes above Dunbar’s number, the limitation on meaningful relationships would result in social networks that increasingly present holes and emerging mesoscale motifs – i.e. more complex interaction graphs – as meaningful relationships become sparser and clustered" \cite{Saavedra2023}. Thus Dunbar's number gives a scale at which there is a build up of societal complexity. In sharing emotional connections, communities of size less than Dunbar's number share a sense of identity. The transition from these communities to acquaintances suggests that these groups that are neither too small and tightly connected nor too large and loosely connected, where the intent to influence is most pertinent. In particular, at the next scale, from Dunbar's number to 1000 (regarded as `small' for online interactions), still supports `purposeful' interactions amongst members, for example `communities of practice', where influence within the common purpose may be desired \cite{Hwang2021}. This is the scale at which we intend our model.

As Dunbar argues, cognitive processing is important in setting the scale for such communities of interaction, and can be considered a dynamic processes that evolves in time \cite{Ditterich2006,Faulkenberry2014}. Stochastic models have been shown to appropriately capture cognitive processes associated with decision making \cite{Voss2009, Voss2013}. For social systems, the importance of stochasticity is critical to capture the heterogeneity of human cognitive processes within an otherwise socially structured dynamical environment. 
Stochastic effects in synchronisation on networks are also well-studied. Both \cite{Bag2007} and \cite{Khoshbakht2008} studied the effects of Gaussian additive noise to the Kuramoto model. The work of \cite{KallZup2013} extended these considerations to additive and multiplicative noise. Here we note that multiplicative noise effectively introduces time-dependent links into the dynamics. However, as inferred above, there are many aspects of interactions in social systems on networks where biased human behaviours recommend heavy tail statistical properties, be they reactive or proactive decision making, enthusiasm or fatigue, or conservation of effort making for frequent weak interactions or infrequent strong interactions. We consider heavy-tailed multiplicative noise as a means of representing influence on a dynamical network; heavy tail distributions for additive noise were first studied in synchronisation by \cite{KallRob2017,RobKall2018}.

In our work stochastically generated dynamical weights are assigned to the influence links, such as in \cite{Brede2016,Kalloniatis2019}, which are comparable to stochastic pinning control formulations seen in recent work \cite{Wu2021,Zhu2024}. We perform numerical experiments to examine transition points in the ability of the influencer to alter the behaviour of the internal network; either to disrupt synchronisation or to drive it to collective frequencies determined by the external node. We examine the dependence of transition points on this node's frequency: too ambitious a frequency fails to influence the system while it retains a synchronised state, and reducing the frequency achieves a state of internal synchronisation, at a collective frequency shifted from the natural mean.  We also look at the analytic approximation for the system close to synchronisation and fracture to understand its behaviour around this limit. We discuss the comparison between the numerical and analytical results to show how the proposed model of influence reflects properties seen in human social systems, and conclude with proposals for future work to further extend the model.

\section{\label{sec:numerical}Characteristic Behaviours for Erd\H{o}s-R\'{e}nyi Random Graphs: Numerical Results}
\subsection{Numerical set-up}
For this model we wish to study the internal synchronisation of the system when it is under such influence, and the ability of the external node to entrain the entire system to its frequency $\Omega$. To study the influence $\theta$ has on the synchronisation of the oscillators, $\varphi_{i}$, we consider the usual Kuramoto order parameter $r$ given by 
\begin{equation}
r = \frac{1}{N} | \sum^{N}_{j=1} e^{i\varphi_{j}} |
\label{eq:order}
\end{equation}
where $0\leq r\leq 1$. Note, $r=1$ when the system is fully synchronised and $r=0$ when there is no synchronisation.

To understand the external node's ability to entrain the network to its frequency, we define a `closeness' ($\Delta$) of the average instantaneous frequency to the external node's frequency, similar to \cite{Kalloniatis2019}. Namely,
\begin{equation}
\Delta = \frac{1}{N}\sum^{N}_{j=1} |\dot{\varphi}_{j}-\Omega|
\label{eq:closeness}
\end{equation}
where $\Delta=0$ indicates the external node has completely pulled the system completely to its frequency $\Omega$. 

Before proceeding to numerical computations we point a modification in this system to a typical property of Kuramoto dynamics. Summing over oscillators $\varphi_i$ in the system we obtain a `sum rule':
\begin{eqnarray} 
\frac{1}{N} \sum_i \dot{\varphi}_{i} = \bar{\omega} &-& \frac{1}{N}\sum_{i\neq j}A_{ij}\left(\frac{1}{k_i} -\frac{1}{k_j}\right)\sin{\left(\varphi_{i}-\varphi_j\right)} \nonumber \\
&-&\frac{\tau}{N} \sum_i
 B_{i} \sin{ \left(\varphi_{i}-\theta\right)} \nonumber \\
\label{eq:sumrule}
\end{eqnarray}
where we note that our choice of exactly symmetric frequencies around zero for computational results have $\bar{\omega}=0$. Thus the noise means the centre of mass of the system is not necessarily zero in the synchronised state; it can be disrupted. Nevertheless, for small choices of shape parameter $\alpha$ we expect typical results for the right hand side of Eq.~(\ref{eq:sumrule}) to vanish when we consider the time-average over an ensemble of Erd\H{o}s-R\'{e}nyi graphs, if the system has decoupled from the external node.

We study this model through numerical solution to equations defined in Eq.~(\ref{eq:model}), using a Runge-Kutta method 4 (10 substeps) and integration time $t_{\text{end}}$ = 300 with time step $\Delta t$ = 0.02. We note that this method uses an adaptive time-step for calculation but returns the solution only for those values specified. The $G_i$ values for the influence links are resampled independently at each output time-step (annealed disorder), representing fluctuating susceptibility to the external influencing node rather than fixed personality differences. Given our focus on randomness in the stochastic elements of the model we use natural frequencies, $\omega_{i}$, of the oscillators that are equally spaced in the interval $[-1,1]$ such that the average natural frequency for any sample instance satisfies exactly $\bar{\omega}=0$. This allows us to study the influence of the external node in regions outside and inside of the interval of natural frequencies, to capture the property of the influencer being variously part of the group versus outside of the group being influenced. Contrastingly, we initialise $\theta$ at 0. In addition we set initial conditions for $\varphi_{i}$ from a uniform distribution on $[-\frac{\pi}{2},\frac{\pi}{2}]$. This is in order to reduce the variance; there can be, even for the standard Kuramoto model, non-zero equilibria that may catch orbits in a frequency (and not phase) synchronised state \cite{Sclosa22}. 

We consider Erd\H{o}s-R\'{e}nyi graphs to capture some social complexity and calculate numerical solutions for an ensemble of graphs ($p=0.3$) with a fully connected external node and calculate ensemble averages using the final 80\% of each time series. We initially baseline the computations for the standard complete graph Kuramoto model where the critical coupling is analytically known and determine for the Erd\H{o}s-R\'{e}nyi case an appropriate 'critical coupling value' in the absence of the influence effect. We seek to see how the influence node impacts on the system when it naturally synchronises. To this end $\sigma=2$ is selected.

We have chosen $N=300$, a large graph but certainly not one that approaches the 'thermodynamic limit' $(N\rightarrow\infty$) typically pursued in complex graph Kuramoto studies. This is consistent with our interest in sizes of systems slightly exceeding the accepted range for Dunbar's number characteristic of more intimate groups or 'purpose-driven' communities \cite{Hwang2021}, as discussed in the introduction. This would represent a challenging but achievable group to influence. However we shall also see that the scaling behaviour of the characteristic transitions in this Influence model are much softer than for the Kuramoto model; the thermodynamic limit is therefore of less significance. We defer further discussion on this to first identify characteristic behaviours. 

\subsection{Order parameter and closeness}
\begin{figure}
    \centering
     \begin{subfigure}[b]{0.49\columnwidth}
         \centering
         \includegraphics[width=\columnwidth]{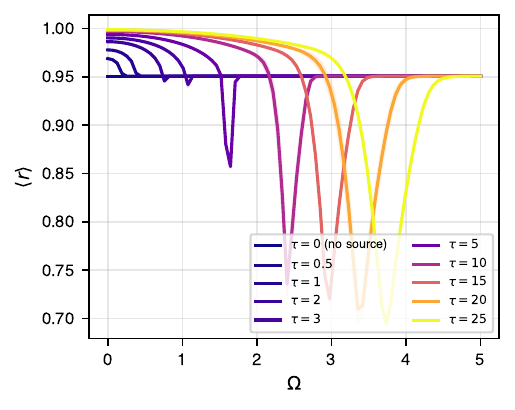}
     \end{subfigure}
     \hfill
     \begin{subfigure}[b]{0.49\columnwidth}
         \centering
         \includegraphics[width=\columnwidth]{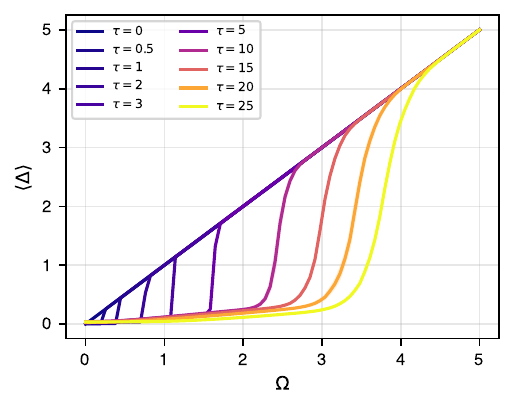}
     \end{subfigure}
    \caption{\textbf{Left: }Time-average order parameter $\langle r\rangle$ vs $\Omega$ (1000 graph realisations, 80 $\Omega$-values in $[0, 5]$). \textbf{Right:} Closeness $\Delta = \frac{1}{N}\sum_i |\langle\dot{\varphi}_i\rangle_t - \Omega|$ vs $\Omega$ (100 graph realisations)on ER($N=300$, $p=0.3$). All graphs with $\tau \in \{0, 0.5, 1, 2, 3, 5, 10, 15, 20, 25\}$, $\sigma = 2$, $\alpha = 1$.}
    \label{fig:regimes}
\end{figure}

We first explore whether mutual synchronisation can be distinguished from entrainment to the external node, through Fig.~\ref{fig:regimes}. Note that $r$ alone cannot distinguish these mechanisms: because $|e^{-i\theta}|=1$, the na\"ive ``source coherence'' $r_s = (1/N)|\sum_j e^{i(\varphi_j - \theta)}| \equiv r$. Thus we must consider both order parameter $r$ and the closeness $\Delta$. In Fig.~\ref{fig:regimes} \textbf{Left} the $\tau = 0$ curve is flat at $\langle r\rangle = 0.9506$, confirming that $\Omega$ only matters through external node coupling. Increasing $\tau$ raises the low-frequency synchronised state toward $\langle r \rangle \approx 1$ and shifts the dip to larger $\Omega$. The minima move from $\langle r \rangle = 0.857$ at $\Omega \approx 1.65$ for $\tau = 5$, to $0.736$ at $\Omega \approx 2.41$ for $\tau = 10$, $0.710$ at $\Omega \approx 3.35$ for $\tau = 20$, and $0.695$ at $\Omega \approx 3.73$ for $\tau = 25$. For weak influence ($\tau \leq  3$), the dip is shallow and the network remains close to its autonomous coherence; by high external node frequencies all curves return to $\langle r \rangle \approx 0.9506$. The shaded bands for the highlighted $\tau$ values remain narrow (maximum standard deviation $\approx 0.022$), so the dip-and-recovery structure is stable across graph realisations. Fig.~\ref{fig:regimes} \textbf{Right} shows that frequency entrainment widens with influence strength. Here we note that as the average $\langle\dot{\varphi}_i\rangle_t$ is taken over the final 80\% of each trajectory, $\Delta \approx 0$ means the internal oscillators are frequency-entrained to the external node. At $\tau=0$ the external node is decoupled and the curve follows the baseline $\Delta\approx\Omega$ ($\Delta(5)=5.00$). Increasing $\tau$ extends the near-zero plateau and shifts the recovery toward the decoupled baseline to larger $\Omega$. Using a 10\%-of-maximum-smoothed-slope criterion, the recovery onset occurs at $\Omega\approx0.57$ for $\tau=2$, $1.46$ for $\tau=5$, $2.15$ for $\tau=10$, $2.91$ for $\tau=20$, and $3.23$ for $\tau=25$; the corresponding steepest recovery points are $\Omega\approx0.82$, $1.71$, $2.47$, $3.42$, and $3.80$. Thus at $\tau=20$ the inflection begins just below $\Omega=3$ and is steepest near $\Omega=3.4$. 

\begin{figure}
    \centering
    \includegraphics[width=\columnwidth]{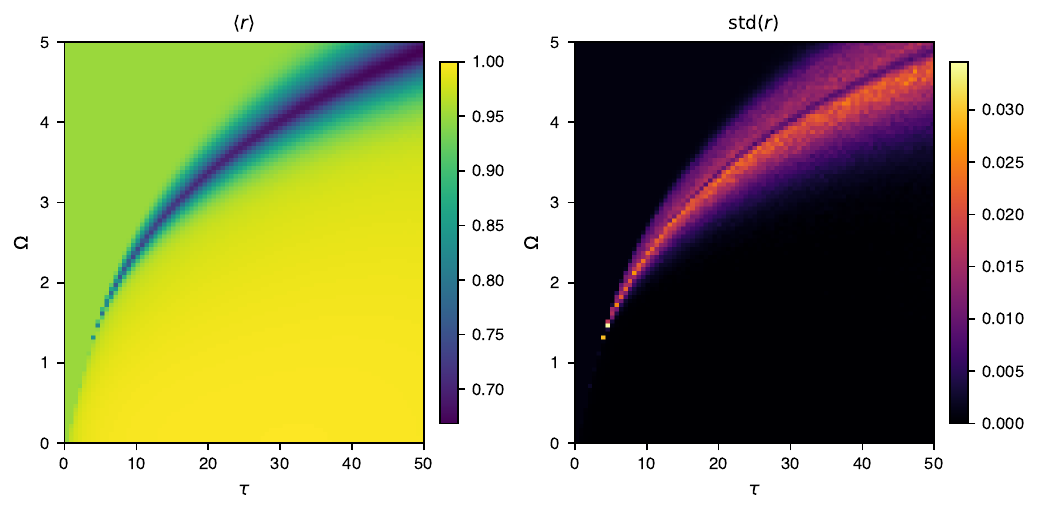}
    \caption{  \textbf{Left:} $\langle r \rangle$ over $\tau \in [0, 50]$ (81 values) and $\Omega \in [0, 5]$ (100 values) on ER($N{=}300$, $p{=}0.3$) graphs ($\sigma=2$, $\alpha=1$; 100 graphs per pixel; 810,000 simulations). \textbf{Right:} $\mathrm{std}(r)$ across graph realisations.}
    \label{fig:order_heatmap}
\end{figure}

Further exploration of model parameters illustrate how the nature of the influence links shape these characteristic behaviours of the model. Fig.~\ref{fig:order_heatmap} shows order parameter $r$ and corresponding standard deviation, examines where exactly the entrained-to-incoherent transition occurs. Here, the autonomous baseline at $\tau=0$ is flat, $\langle r\rangle=0.9508$. Fig.~\ref{fig:order_heatmap} \textbf{Left} shows that for $\tau>0$, the synchronised regime (yellow, $\langle r\rangle\approx1$) occupies the low-$\Omega$ side of the diagram, while the system returns to the autonomous baseline at sufficiently high $\Omega$. Between them is a curved band of reduced coherence. Its centre moves right as influence strengthens: the order-parameter minima occur near $(\tau,\Omega)=(5,1.62)$ with $\langle r\rangle=0.793$, $(10,2.42)$ with $\langle r\rangle=0.740$, $(20,3.38)$ with $\langle r\rangle=0.703$, and $(50,4.90)$ with $\langle r\rangle=0.669$. Thus $\tau=50$ already pushes the transition to the upper edge of the $\Omega\le5$ window; larger $\tau$ would mainly move the same boundary beyond the plotted range. Fig.~\ref{fig:order_heatmap} \textbf{Right} shows that the variability band follows the same curved boundary, peaking at $\mathrm{std}(r)\approx0.035$ near $(\tau,\Omega)\approx(4.4,1.46)$ and remaining much smaller than in the earlier $N=40$ phase diagram. Away from the boundary, graph-to-graph variability is low because realisations agree on either the synchronised high-coherence state or the decoupled autonomous baseline. We note here the band of low variance in \textbf{Right} aligns to the band of lowest order parameter $\langle r \rangle$ in \textbf{Left}. This indicates that this system is showing uniformly low values in this arc, across graph instances. That is, the system is losing coherence in a consistent way.

\begin{figure}
    \centering
    \includegraphics[width=\columnwidth]{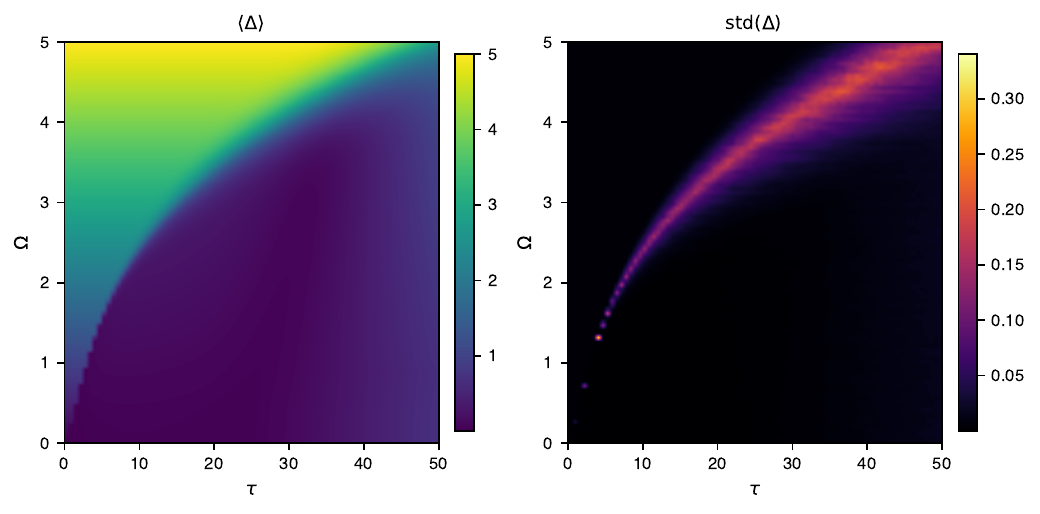}
    \caption{ \textbf{Left:} Closeness $\langle\Delta\rangle = \frac{1}{N} \sum_i |\langle\dot{\varphi}_i\rangle_t - \Omega|$ over $\tau \in [0, 50]$ (81 values) and $\Omega \in [0, 5]$ (100 values), on ER($N{=}300$, $p{=}0.3$) graphs ($\sigma=2$, $\alpha=1$, 100 graphs per pixel). 
    \textbf{Right:} $\mathrm{std}(\Delta)$ across graph realisations.}
    \label{fig:closeness_heatmap}
\end{figure}

Looking at the closeness, Fig.~\ref{fig:closeness_heatmap} explores if the frequency deviation shows the same boundary structure as the order parameter. Noting at $\tau=0$, the external node is ignored and $\langle\Delta\rangle\approx\Omega$. Fig.~\ref{fig:closeness_heatmap} \textbf{Left} shows increasing $\tau$ opens a low-$\Delta$ entrained region whose recovery boundary tracks the order-parameter band in Fig.~\ref{fig:order_heatmap}: the heuristic recovery onset moves from $\Omega\approx1.46$ at $\tau=5$ to $2.17$ at $\tau=10$, $2.93$ at $\tau=20$, and $3.94$ at $\tau=40$. For very strong influence ($\tau\gtrsim40$), low-frequency $\Delta$ rises again even though $\langle r\rangle$ remains high, indicating a coherent but over-stretched response rather than simple frequency locking. Fig.~\ref{fig:closeness_heatmap} \textbf{Right} shows that the variability again concentrates along the entrainment boundary, with a maximum $\mathrm{std}(\Delta)\approx0.34$ near $(\tau,\Omega)\approx(3.75,1.31)$ and smaller but persistent bands at larger $\tau$. Thus $\Delta$ separates frequency entrainment from mutual coherence: it agrees with the order-parameter boundary but also reveals where strong influence keeps oscillators coherent without making their mean frequencies equal to the external node. The closeness metric provides a complementary view of the same phase structure, directly measuring frequency entrainment rather than phase coherence. 

\subsection{Shape of oscillator configurations under influence}
\begin{figure*}
    \centering
    \includegraphics{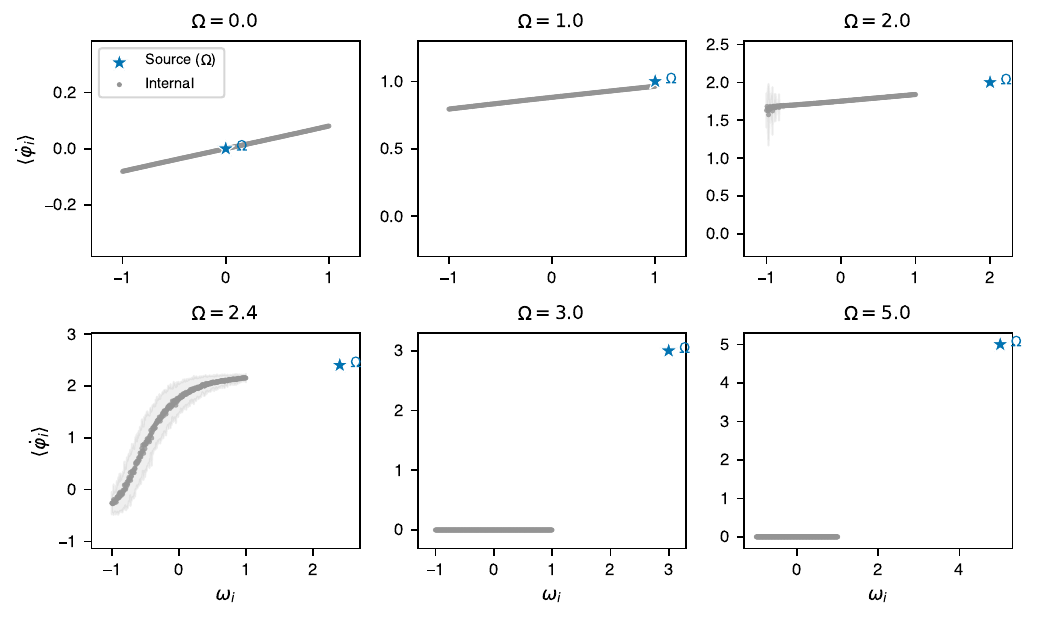}
    \caption{Ensemble-averaged effective frequency $\langle\dot{\varphi}_i\rangle_t$ (grey dots) for each internal oscillator, plotted against its natural frequency $\omega_i$, across 100 ER($N{=}300$, $p{=}0.3$) graphs, with $\sigma=2$, $\tau=10$, $\alpha=1$; $\Omega \in \{0, 1, 2, 2.4, 3, 5\}$. The shaded band shows $\pm 1\sigma$ across graph realisations; the blue star marks the external node at frequency $\Omega$. }
    \label{fig:stretching}
\end{figure*}

We now turn to examining which oscillators get pulled toward the external frequency. We look to identify patterns of stretching and whether fracture occurs within the system due to the influencing node. To do so we consider time-averaged frequencies. Fig.~\ref{fig:stretching} shows that at $\Omega=0$ the internal population remains close to the stationary external node, with only a narrow spread around zero. At $\Omega=1$ the entire oscillator band is pulled upward and remains tightly organised ($\langle\dot{\varphi}_i\rangle_t \approx 0.8\text{--}1.0$). At $\Omega=2$ the cloud is stretched toward the external node and the variance grows, indicating the onset of graph-dependent partial entrainment. The most extremely stretched case is $\Omega=2.4$, matching the minimum of the $N=300$, $\tau=10$ order-parameter curve in Fig.~\ref{fig:regimes} \textbf{Left} oscillator frequencies spread from slightly negative values up toward the external node, with the largest inter-graph uncertainty. By $\Omega=3$ and $5$ the external node has visibly decoupled, and the internal oscillators collapse back near their autonomous mean frequency $\langle\dot{\varphi}_i\rangle_t \approx 0$. 

\subsection{Scaling behaviour}
\begin{figure}
    \centering
    \includegraphics[width=\linewidth]{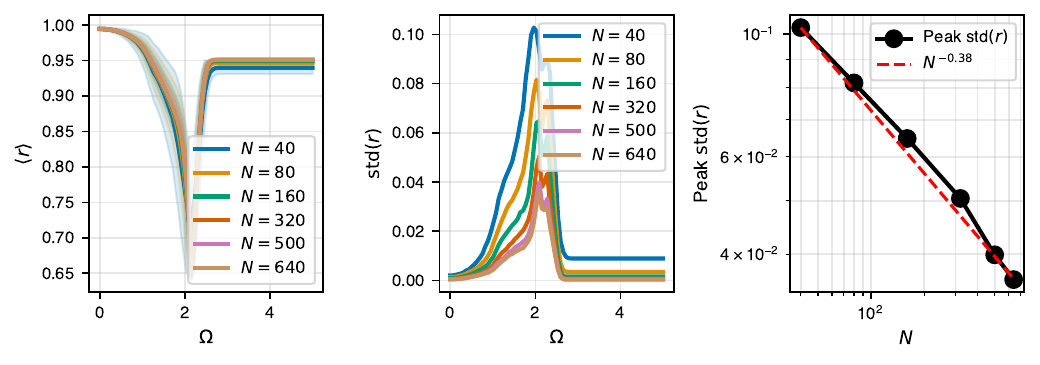}
    \caption{
    \textbf{Left:} $\langle r \rangle$ vs $\Omega$ for $N \in \{40, 80, 160, 320, 500, 640\}$ at $\tau=10$, $\alpha=1$; 500--1000 ER($p{=}0.3$) graphs per $N$, 80 $\Omega$-values in $[0, 5]$. Shaded bands: $\pm 1$ std.\ dev. 
    \textbf{Centre:} $\mathrm{std}(r)$ vs $\Omega$. 
    \textbf{Right:} Peak $\mathrm{std}(r)$ vs $N$ on log--log axes. }
    \label{fig:scaling_orderparameter}
\end{figure}

Having now identified the key behaviours of the model we examine their finite size scaling to test whether they are critical in the same sense as that for the Kuramoto order parameter. Fig.~\ref{fig:scaling_orderparameter} explores if the entrainment transition sharpens with system size. Here Fig.~\ref{fig:scaling_orderparameter} \textbf{Left} illustrates that all sizes show the dip-and-recovery pattern with the dip near $\Omega \approx 2.15$. The dip deepens modestly from $\langle r \rangle \approx 0.73$ ($N=40$) to $\langle r \rangle \approx 0.67$ ($N=640$). The transition steepens markedly with $N$: the $N=640$ curve descends almost vertically near $\Omega \approx 1.8$ and recovers sharply by $\Omega \approx 2.5$, while the $N=40$ curve shows a broader, more gradual dip. The decoupled baseline converges to $\langle r \rangle \approx 0.95$ for all $N \geq 160$. Curves for $N = 500$ and $N = 640$ are nearly indistinguishable, indicating convergence. Fig.~\ref{fig:scaling_orderparameter} \textbf{Centre} shows the fluctuations peak in the stretched  regime, decreasing from $\sim\!0.10$ ($N=40$) to $\sim\!0.036$ ($N=640$). The peak narrows with $N$ as the transition sharpens. Finally, Fig.~\ref{fig:scaling_orderparameter} \textbf{Right} shows the fit gives $\mathrm{std}(r) \propto N^{-0.38}$ (95\% CI: $[-0.41, -0.35]$). Note the exponents are shallower than $N^{-1/2}$, reflecting intra-graph correlations in the synchronised state.
We see therefore that in the presence of the influence effect 
the `critical' transition is much softer than that for the Kuramoto model; there is no abrupt system wide change when the system couples to the driver.

\subsection{Deterministic limit}
\begin{figure}
    \centering
    \includegraphics[width=1\columnwidth]{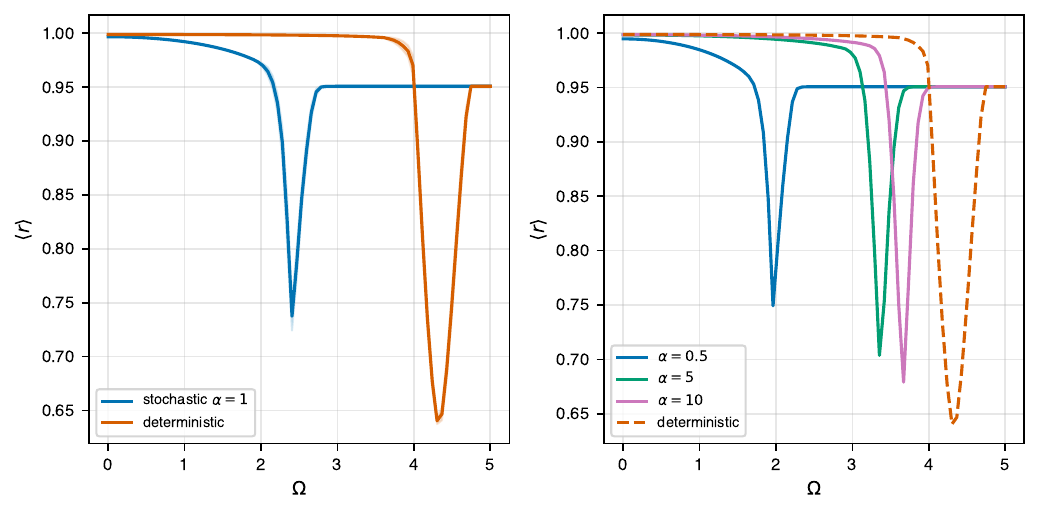}
    \caption{$\langle r \rangle$ vs $\Omega$ on ER($N=300$, $p=0.3$) graphs, $\sigma = 2$, $\tau = 10$, 100 graph realisations, 80 $\Omega$-values in $[0, 5]$ \textbf{Left:} Comparing stochastic links ($G_i \sim \Gamma(1, 1)$, blue) with deterministic links ($Gi \equiv 1$, vermillion), \textbf{Right:} Comparing stochastic links ($G_i \sim \Gamma(\alpha, 1)$: $\alpha\in\{0.5,5,10\}$) with deterministic links ($Gi \equiv 1$)}
    \label{fig:deterministicR}
\end{figure}

We consider now whether the randomness in the influence links matters, particularly how it affects the characteristic behaviour of $\langle r\rangle$. Looking at Fig.~\ref{fig:deterministicR}, we see that both curves in \textbf{Left:} show high coherence for low $\Omega \lesssim 1$, then decline until they reach a minimum and return steadily to the autonomous baseline. For stochastic influence links $\alpha = 1$, we see this drops earlier but less deeply than the deterministic limit: the stochastic curve reaches its minimum $\langle r\rangle \approx 0.738$ at $\Omega\approx 2.41$, while the deterministic curve reaches a deeper minimum $\langle r \rangle \approx 0.640$ at $\Omega \approx 4.30$. Thus indicating that this regime is both delayed and more severe without stochastic variability. In Fig.~\ref{fig:deterministicR} \textbf{Right:} we note that increasing $\alpha$ shifts the dip in $\langle r \rangle$ toward the deterministic case and also deepens it: the minima move from $\Omega \approx 1.96$ at $\alpha = 0.5$, to 3.35 at $\alpha = 5$, and 3.67 at $\alpha = 10$. Thus variability in the influence weights advances the onset of coherence loss but softens the de-cohered state. Noting that for increasing $\alpha$, the Gamma distribution concentrates near its mean, the stochastic and deterministic cases converge. Thus confirming that the stochastic–deterministic gap is driven by the variance of $G_i$. 
\begin{figure}
    \centering
    \includegraphics[width=1\columnwidth]{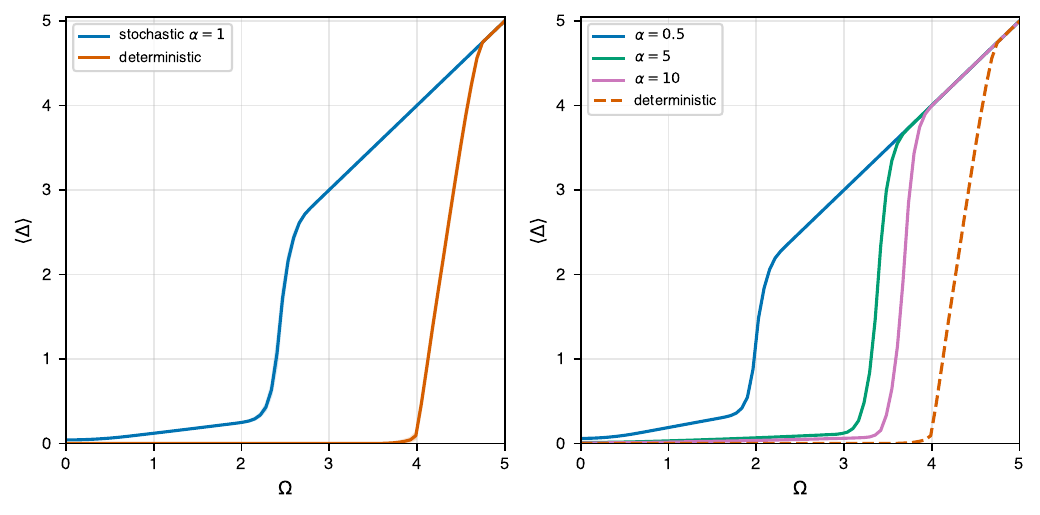}
    \caption{$\langle\Delta\rangle$ vs $\Omega$ on ER($N=300$, $p=0.3$) graphs, $\sigma = 2$, $\tau = 10$, 100 graph realisations, 80 $\Omega$-values in $[0, 5]$ \textbf{Left:} Comparing stochastic links ($G_i \sim \Gamma(1, 1)$, blue) with deterministic links ($Gi \equiv 1$, vermillion), \textbf{Right:} Comparing stochastic links ($G_i \sim \Gamma(\alpha, 1)$: $\alpha\in\{0.5,5,10\}$) with deterministic links ($Gi \equiv 1$)}
    \label{fig:deterministicDelta}
\end{figure}
The stochastic–deterministic separation is also visible in the closeness, as in Fig.~\ref{fig:deterministicDelta}. \textbf{Left} illustrates that stochastic $\alpha = 1$ loses source-frequency entrainment much earlier than the deterministic limit: the steepest rise in $\langle\Delta\rangle$ occurs near $\Omega \approx 2.47$ for $\alpha = 1$, but near $\Omega \approx 4.24$ for deterministic influence. At $\Omega = 3$, the stochastic curve has already returned to the decoupled value $\langle\Delta\rangle \approx3$, while the deterministic curve remains essentially entrained with $\langle\Delta\rangle \approx 0$. Fig.~\ref{fig:deterministicDelta} \textbf{Right} shows that increasing $\alpha$ delays the loss of entrainment, with steepest-rise locations moving from $\Omega \approx 2.03$ at $\alpha = 0.5$, to 3.35 at $\alpha = 5$, and 3.67 at $\alpha = 10$, approaching the deterministic curve at $\Omega \approx 4.24$. Thus the Delta diagnostic confirms the order-parameter picture in Fig.~\ref{fig:deterministicR}: influence-weight variability advances the transition out of source locking. 

\subsection{Summary of behaviours}
To summarise the characteristic behaviours, we see that this model of influence demonstrates regimes in which the external node is able to entrain the whole system while it maintains synchronisation, stretch the system such that the order parameter decreases, and ultimately decouple completely from the external node and return to it's autonomous, internal synchronisation. The model shows a smooth, not sharp, phase transition and thus is of relevance to systems of a size relevant to social community systems. 

\section{\label{sec:analytics}Analytic Results via Linear Approximation}
\subsection{Ansatz for clustering and synchronisation}
We undertake a stability analysis, to understand the phenomenon of entrainment and clustering in this system analytically and whether the external node can cause fracture of the internal system into two clusters. We first assume a system that interpolates between coherence and fracture, such that we can consider two clusters of oscillators $\Phi_{1}$ and $\Phi_{2}$ with 
\begin{align} 
\Phi_{1} &= \theta + \delta \nonumber \\
\Phi_{2} &= \Phi_{1} + \epsilon 
\label{eq:assume}
\end{align}
where $\Phi_{1}$ and $\Phi_{2}$ have mean frequency $\bar{\omega}_{1}$ and $\bar{\omega}_{2}$ respectively, their phase separation is $\epsilon$, and we define $\delta$ as the phase separation between $\theta$ and $\Phi_{1}$. In doing so, we are able to analytically study the system in the regime of both entrainment and potential fracture, where $\epsilon$ is representative of the degree of separation of clusters within the internal system and $\delta$ is connected to the capture by the external node. 
We write Eq.~(\ref{eq:model}) as
\begin{align} 
\dot{\theta} &= \Omega \nonumber \\ 
\dot{\Phi}_{1} &= \bar{\omega}_{1} - \sigma \sin{ \left(\Phi_{1}-\Phi_{2}\right)} - \tau \bar{B}_{1} \sin{ \left(\Phi_{1}-\theta\right)}\label{eq:frac1}\\
\dot{\Phi}_{2} &= \bar{\omega}_{2} - \sigma \sin{ \left(\Phi_{2}-\Phi_{1}\right)} - \tau \bar{B}_{2} \sin{ \left(\Phi_{2}-\theta\right)}\label{eq:frac2}
\end{align}
Here $\delta$ and $\epsilon$ are fluctuations; these will not necessarily be taken from the outset as `small'. Rewriting Eq.~(\ref{eq:frac1}) and Eq.~(\ref{eq:frac2}) we have
\begin{align} 
\dot{\delta} = &\bar{\omega}_{1} - \Omega+ \sigma \sin(\epsilon) - \tau  \bar{B}_{1} \sin{ \left(\delta\right)}\label{eq:fracdelta} \\
\dot{\epsilon} = &\left(\bar{\omega}_{2} -\bar{\omega}_{1}\right) - 2\sigma \sin(\epsilon) \nonumber\\
&- \tau\left[ \bar{B}_{2} (\sin \left(\delta\right)\cos(\epsilon)+\cos(\delta)\sin(\epsilon)) - \bar{B}_{1} \sin \left(\delta\right)\right]\label{eq:fracepsilon}
\end{align}
We now include the explicit form of $\bar{B}_1$ and $\bar{B}_2$, defined by $\bar{B}_{a}=\frac{1}{N_a}\sum_{i\in a}\frac{G_{i}}{(\varphi_i-\theta)^{2}+1}$, where $a\in\{1,2\}$ and $G_{i}\sim\Gamma(\alpha,\frac{1}{\alpha})$. Thus, we use the expressions
\begin{align}
\bar{B}_{1}&=\frac{\bar{G}_{1}}{\delta^{2}+1} \\
\bar{B}_{2}&=\frac{\bar{G}_{2}}{(\delta+\epsilon)^{2}+1},
\end{align}
where $\bar{G}_{1,2}$ represent the sum over the noise terms across clusters 1 and 2 normalised by the size of the clusters.
Assuming $\Phi_{1}$ and $\Phi_{2}$ are nearly phase-synchronised ($\Phi_{1}-\Phi_{2}=\epsilon$ is small), and close to capture ($\Phi_{1} -\theta = \delta$ is small) we are able to write the system as follows
\begin{align} 
\dot{\delta} = &\bar{\omega}_{1} - \Omega+ \sigma \epsilon - \tau  \bar{G}_{1} \delta+ O(\delta^{2})+O(\epsilon^{2})\label{eq:lindelta} \\
\dot{\epsilon} = &\left(\bar{\omega}_{2} -\bar{\omega}_{1}\right) - \left(2\sigma +\tau\bar{G}_2\right)\epsilon \nonumber\\
&- \tau\left( \bar{G}_{2} - \bar{G}_{1}\right)\delta+ O(\delta^{2})+O(\epsilon^{2})\label{eq:linepsilon}
\end{align}
Thus we obtain the linear system as
\begin{equation}
\begin{gathered}
 \begin{bmatrix} \dot{\epsilon} \\ \dot{\delta} \end{bmatrix}
 =\begin{bmatrix} \bar{\omega}_2-\bar{\omega}_1 \\ \bar{\omega}_1-\Omega \end{bmatrix}-
  \begin{bmatrix}
   2\sigma+\tau \bar{G}_2  &
   \tau (\bar{G}_2-\bar{G}_1)  \\
   -\sigma  &
    \tau \bar{G}_1  
   \end{bmatrix}
   \begin{bmatrix} \epsilon \\ \delta \end{bmatrix}
   \label{eq:matrix}
\end{gathered}
\end{equation}
Note that time-dependence remains in this system through the noise $\bar{G}_{1,2}$, which are sampled continuously in time. The full solutions are given in Appendix C.

\subsection{Stability analysis of the cluster ansatz}
We proceed in order to understand the qualitative behaviours of the system with regard to stability and equilibrium behaviours based on the Lyapunov eigenvalue spectrum of the matrix in Eq.(\ref{eq:matrix}). The eigenvalues are given by
\begin{widetext}
\begin{equation}
\lambda_{\pm}=\frac{1}{2}\left( 2\sigma+\tau(\bar{G}_1+\bar{G}_2)\pm\sqrt{\left(2\sigma+\tau(\bar{G}_1+\bar{G}_2)\right)^2-4\tau\left((\bar{G}_2+\bar{G}_1)\sigma+\bar{G}_1\bar{G}_2\tau\right)} \right)
\label{eq:lambda}
\end{equation}
\end{widetext}
This can be written as
\begin{equation}
\lambda_{\pm}=\frac{1}{2}\left( 2\sigma+\tau(\bar{G}_1+\bar{G}_2)\right) \left( 1 \pm\sqrt{1- 4\tau H} \right)
\label{eq:lambdaH}
\end{equation}
with 
\begin{equation}
    H \equiv {\left( (\bar{G}_1 + \bar{G}_2) \sigma + \tau \bar{G}_1\bar{G}_2 \right) \over \left( 2\sigma+\tau(\bar{G}_1+\bar{G}_2) \right)}.
\end{equation}
The form of Eq.(\ref{eq:lambdaH}) reveals that $\mathbb{R}(\lambda)>0$ thus this system always has stable solutions (noting our sign convention). 
The eigenvalues remain time-dependent through the noise terms thus they appear through time-integrals in the solutions, as explained in Appendix C. We distinguish three cases: $\tau=0$, $0<H<\frac{1}{4\tau}$ -- namely $H$ is small -- and $H \geq 1/4\tau$.
For $\tau=0$ we obtain $\lambda=\{0,2\sigma\}$ where the zero mode corresponds to the {\it conserved} centre of mass frame, and the non-zero mode the stable equilibrium decaying back to the Kuramoto synchronised state (discussed in Appendix C).

Through the Gamma distribution, $G_{1,2}$ will more frequently draw small values less than one given our choice of $\alpha=1$ as in Section \ref{sec:numerical}; the heavy tail in $\Gamma$ means that, less frequently, large values may arise. Thus, {\it typically} $\bar{G}_1\bar{G}_2$ will be vanishingly small. Therefore
\begin{eqnarray}
    H &\approx & \sigma { ( \bar{G}_1 + \bar{G}_2 ) \over 
     \left( 2\sigma+\tau(\bar{G}_1+\bar{G}_2) \right)} \nonumber \\
     &\approx & \frac{1}{2} { ( \bar{G}_1 + \bar{G}_2 ) \over 
     \left( 1 +\tau(\bar{G}_1+\bar{G}_2)/2\sigma \right)} \nonumber \\
     &\approx & \frac{1}{2} ( \bar{G}_1 + \bar{G}_2 ) 
\label{eq:Happrox}
\end{eqnarray}
to leading order in $G_{1,2}$ for the typically small noise samples. 
Significantly this is $\sigma$ independent, so that $H$ will typically be small {\it regardless} of the underlying Kuramoto system parameters. 
In this case the eigenvalues are
\begin{equation}
\lambda = \left(2\sigma + \tau (\bar{G}_1 + \bar{G}_2 )\right) \{H, 1-H\}
\label{lambdaHweak}
\end{equation}
where evidently $H>0$. We see here that the $\sigma$ dependence is purely additive, leading to an overall multiplicative coefficient in the damping in the full solution (seen in Appendix C). 
Indeed, for two approximately equal fragments where both scale similarly in $N$, the sampling of the noise will be similar so that $\bar{G}_1\approx \bar{G}_2$. The form of the solution in Appendix C
shows it collapsing to the standard Kuramoto form reflective of
the centre of mass and relative motion, but with enhanced suppression through the positive eigenvalues of Eq.(\ref{lambdaHweak}). 

Less frequently there will be a draw of the Gamma noise that will be large, so $H>\frac{1}{4\tau}$ and the system will have imaginary components, $\mathbb{I}(\lambda)\neq0$.
This means that there will be noise instances where, although stable, there will be time-dependent oscillations. In other words, 
when such a kick occurs one or both of $\dot{\epsilon}$ and $\dot{\delta}$
will be non-zero in the presence of the influencer over a period of time. The larger $\tau$ is, within the bounds of validity of the approximation giving Eq.(\ref{eq:Happrox}), the smaller $\frac{1}{4}\tau$ and therefore the less the tolerance of the system for kicks of even smaller size. 

If $\tau$ is arbitrarily large in relation to $\sigma$ and the scale of noise we may approximate directly in the original form of the eigenvalue
Eq.(\ref{eq:lambda}) by retaining terms of order $\tau^2$ under the square root:
\begin{eqnarray}
\lambda_{\pm} & \approx & \frac{1}{2}\left( \tau (\bar{G}_1 + \bar{G}_2) 
\pm \sqrt{\tau^2 (\bar{G}_1 + \bar{G}_2)^2 - 4 \tau^2 \bar{G}_1 \bar{G}_2} \right) \nonumber \\
&=& \frac{\tau}{2}\left( (\bar{G}_1 + \bar{G}_2) 
\pm | \bar{G}_1 - \bar{G}_2 | \right)
\end{eqnarray}
thus $\lambda=\{\tau\bar{G}_1,\tau \bar{G}_2 \}$ for, say, $\bar{G}_1 >\bar{G}_2$ .
There is no imaginary part. This state is therefore
relatively {\it more stable} than the smaller $\tau$ case.

Overall, comparing Eq.(\ref{lambdaHweak}) with $H\neq 0$ to $\{ 0,2\sigma \}$ we conclude that the typical behaviour of the system in the presence of the influencer is that fluctuations in $\epsilon$ (clusters within the population) and $\delta$ (deviations between the influencer and the population) are {\it more suppressed} than for the standard Kuramoto system left to itself. With draws of the noise that are small, there is also no zero mode and therefore no constant direction in the fluctuations - consistent with the lack of conservation of the centre of mass of the system. {\it Nevertheless, entrainment of the system by the influencer is Lyapunov stable.} Moreover, any mode where fracture between $N_1$ and $N_2$ might occur is quickly suppressed; fluctuations in the inter-cluster phase difference $\epsilon$ will rapidly decay. 
These statements are localised in time until the infrequent instant where there is a large sample from the noise, a kick. Then an imaginary part generates an oscillation where potentially the `smallness' of the solution is violated, the system leaves the basin of attraction until recapture, with continuing small nudges then being damped again.
Contrastingly, increasingly large $\tau$ generates greater sensitivity to kicks to smaller and smaller scales up until $\tau$ compensates for
the noise of order $\bar{G}^2$. Then relative stability is recovered.

\subsubsection{Stochastic vs deterministic influence links}
For the deterministic case, which corresponds to $\bar{G}_1=\bar{G}_2=1$ the eigenvalues are simply $\lambda=\{ \tau, 2\sigma + \tau \}$ (thus stable). But as shown in Appendix C, the coefficients in the solutions to the linear system reduce to the ordinary Kuramoto case; this more as a consequence of their dependence on the noise difference, $\delta G \equiv \bar{G}_1-\bar{G}_2$ =0. However, for the stochastic system for any time instant because of the heavy tail the noise sample in one cluster may be very different from the other cluster; $\delta G$ may be large. Taking leading terms in $\delta G$ in the coefficients $\xi_D$ sees the $G$ dependence cancel in the numerator but not in the denominator, which gives
\begin{equation}
\xi_{D_\pm} = \pm { 2\sigma (\Omega-\bar{\omega}_2) \over {2\tau \delta G}}.
\end{equation}
Thus the coefficients in the equilibrium solutions for the stochastic case are suppressed by the noise terms compared to the deterministic case.

Qualitatively, these behaviours are consistent with the overall numerical results. Specifically, in the heatmaps of Figs.\ref{fig:order_heatmap} and \ref{fig:closeness_heatmap} we see that for any fixed $\Omega$ the order parameter and closeness degrade for increasing $\tau$. There is then a recovery in line with a strong $\tau$ limit. The suppression of fracture configurations is also consistent with the softness of the entrainment transition from the scaling analysis. Finally, stability analysis shows that the equilibrium solutions for the stochastic case will be more suppressed than for the deterministic case.

\subsection{Order parameter and closeness}
Returning now to our definition of the order parameter and closeness, as given in Eq.~(\ref{eq:order}) and Eq.~(\ref{eq:closeness}) respectively. We can write these in terms of $\epsilon$ and $\delta$. The order parameter is written as 
\begin{equation}
r \approx\frac{1}{N_{1}+N_{2}} \sqrt{N_{1}^{2}+N_{2}^{2}+2N_1N_2\cos(\epsilon)}
\end{equation}
We note that when $\epsilon=0$, and the internal system is fully synchronised then $r=1$ as expected. Here, $\epsilon$ is given by combinations of the general solutions $X_{D_{\mp}}$ as given in Appendix \ref{app:solution}. We note also that as $\epsilon$ increases ($\dot{\epsilon}\neq0$), the system starts to stretch and becomes less synchronised, thus the order parameter decreases as expected (as represented in Fig.~\ref{fig:regimes} \textbf{Left}). While the expression for $r$ only depends explicitly on $\epsilon$, we note that the interconnected dynamics between $\epsilon$ and $\delta$ has some interplay with regards to the internal synchronisation. 

Considering the closeness now, we have
\begin{equation}
\Delta \leq \frac{1}{N}\left(N_{1}|\dot{\delta}|+N_{2}|\dot{\delta}+\dot{\epsilon}|\right)
\end{equation}
where the RHS of Eq.~\ref{eq:matrix} and general solutions in Appendix \ref{app:solution} can be inserted. We consider $\Delta$ in multiple regimes in order to understand the behaviour of the system that is represented in the numerical results in Fig.~\ref{fig:regimes} \textbf{Right}. If the system is both entrained to the external node and remains internally synchronised then $\Delta\rightarrow0$ as $\dot{\epsilon}\rightarrow0$ and $\dot{\delta}\rightarrow0$. In this state, the external node may be seen as \textit{part of the system}, and thus synchronises with the other nodes. This is seen in the plateau at $\Delta\approx0$ in Fig.~\ref{fig:regimes} \textbf{Right}. There is a point however at which the external node may be seen as \textit{influencing the system} towards its own behaviour, namely the system no longer synchronises to $\bar{\omega}=0$. We note that for a system that is phase-locked to the external node, such that $\dot{\delta}\rightarrow0$, the closeness is bound by $\Delta\leq\frac{N_{2}}{N_1+N_2}|\dot{\epsilon}|$. Alternatively, when the system is synchronised internally and not synchronised to the external node, that is where $\dot{\epsilon}\rightarrow0$, the closeness is bound by $\Delta\leq|\dot{\delta}|$. As $\Omega$ becomes larger the system ignores the external influencer and resynchronises to $\bar{\omega}=0$. Specifically, in closeness we see that the inflection away from the plateau in Fig.~\ref{fig:regimes} \textbf{Right} corresponds to $\dot{\epsilon}\neq0$ and thus the system here is stretched, as described by $r$ anove. Analytically, we see that as $\Omega$ becomes large $\Delta\leq \Omega$, from Eq.~(\ref{eq:lindelta}), and see computationally that $\Delta=\Omega$ for large $\Omega$. Thus the analytical approximation leads to three distinct regimes in $\Delta$: captured, stretched and decoupled, consistent with behaviours seen in Fig.~\ref{fig:regimes}.

As described above, the tendency for one or both of $\dot{\epsilon}$ and $\dot{\delta}$ to be non-zero creates a highly dynamic influence on the internal system. We see that $\dot{\epsilon}=\dot{\delta}=0$ represents \textit{entrainment}; $\dot{\epsilon}\neq0$ and $\dot{\delta}=0$ represents a \textit{stretched} system; and, $\dot{\epsilon}=0$ and $\dot{\delta}\neq0$ represents a \textit{decoupled} system. These are indicative now of the concepts for the behavioural regimes for this model: a system that is \textit{entrained}, \textit{stretched}, or \textit{decoupled} from the influencing node. 

\section{\label{sec:discussion}Insights and Discussion}
Here we consider the characteristic behaviours of the system with respect to the influencing node. These show that the system transitions through a captured regime, where it remains synchronised and entrained to the external node. In the captured regime (small $\Omega$), $r$ is high and $\Delta$ is low because oscillators are mutually synchronised \emph{and} entrained to the external node. Contrastingly, in the decoupled regime (large $\Omega$), $r$ returns to the autonomous baseline while $\Delta$ increases to $\Omega$: oscillators are synchronised to each other but not to the external node. Between these is a regime where $\langle r\rangle$ dips, and $\Delta$ begins to increase from the plateau, where the system is stretched. One can consider it to be maximally stretched at the lowest point in $\langle r\rangle$, where it is typically partly decoupled. The intrinsically dynamic nature of this model for influence is captured both numerical and analytically in this paper. Consolidating the results in section \ref{sec:numerical} and \ref{sec:analytics}, we see a three-regime interpretation for this model: \textit{captured}, \textit{stretched} and \textit{decoupled}.

The relationship between the external node on the collective system can be understood as {\it influence} through the collation of the above results, both analytic and numeric. Overall we observe a system that internally synchronises, but is dynamically influenced by the external node with different behaviours depending on whether the driving frequency $\Omega$ places the driver within the group (as defined by the natural frequency range of the oscillators), at or close to its edge, or well outside the group. The system synchronises to the external node's frequency $\Omega\neq\bar{\omega}$ for small values of $\Omega$, namely when the external node is \textit{part of the group}, well within the range of the system's natural frequencies. As $\Omega$ moves to and passes the edge of the group, it influences the system such that it develops into a \textit{stretched} state and is pulled toward $\Omega$. Here we consider a stretched state to be when the nodes remain synchronised but are more dispersed and more loosely bound than they otherwise would be. This represents the least synchronised state for the system, where the tension between internal synchronisation and capture to the external node is at its strongest, as is evident in the order-parameter at its lowest value r (see Fig.~\ref{fig:regimes}). Depending on the parameters of the model, the internal system can be dropped by the external node and quickly resynchronise to its natural mean frequency, $\bar{\omega}=0$ (subject to the average noise in Eq.~(\ref{eq:sumrule})), or it will capture for some time before ultimately being dropped. In the case of fracturing, all nodes are slowly dropped over time and resynchronise to the natural mean. This is consistent with the analytic findings in Section $\ref{sec:analytics}$

There is therefore a niche parameter regime in which the system will be in its maximally stretched state, in which it is pulled towards the external node before drifting back towards its natural mean frequency $\bar{\omega}$ as it decouples. {\it This is the region of maximum influence}, and is identified by the minimum in the order parameter and a corresponding increase from the plateau of the closeness in Fig.~\ref{fig:regimes}. We see that the nodes that experience its influence the strongest, are those closest in natural frequency to the external node. The soft transition shows that that the system is flexible under the influence of the external node, and thus it stretches. Our results show that the external node does not cause further fracture to an already synchronised system. 

As described above, the model also exhibits \textit{limitations} on the influence that the external node can have on the internal system. As the influence node drops the system after a finite time, there is a window of influence for which it can ultimately \textit{capture} the system, which is dependent on parameter choices of the model that shape the influence links. We note also that this window is larger for deterministic links but the constant interaction with the external node means that the system does not re-synchronise after it is dropped. Comparatively, this window is smaller for the stochastic links but due to the comparatively light and intermittent interactions of the stochastic links, the system is able to re-synchronise.

Poignantly, this model demonstrates that a model of influence can both capture the aspect of \textit{identity} and \textit{autonomy} as described in \cite{West2025}. The different regimes of behaviour and be seen as dependent on whether the external node is inside the population or not, and the system being correspondingly \textit{captured}, \textit{stretched} or \textit{decoupled} illustrates the subtlety in capturing elements of how the influencer speaks to identity within a group. In particular, the stretching of the internal system is representative of a destabilised sense of cohesion or deepening divisions within a population, as discussed in \cite{West2025,West2023}. The profile of the influence links ultimately governs the influence on the system, particularly if they are stochastic in nature rather than deterministic. Specifically, the stochastic links allow the system to maintain a higher level of internal synchronisation and thus {\it autonomy of the group to naturally interact within itself} across a range of values for $\Omega$. This is most prominently seen for $\alpha=1$ where the internal synchronisation is able to be maintained while the system is drawn towards the external node; very weak link strengths are more frequently activated, which may be described as a `light, deft touch' as most desirable for effective social influence. Larger values of $\alpha$ (tending towards the deterministic case) with more commonly non-zero link strengths activated by the external node, destroy the ability of the system to retain high levels of internal synchronisation. This reinforces the importance of stochastic `soft' links in an influence model, specifically in order to maintain the autonomy of the social system. This model of influence speaks to the nuanced balance between guiding and shaping the synchronisation of the internal system by drawing it to the external node, while the internal system maintains it's natural ability to synchronise.

\section{Conclusions}
Summarising the results of this paper, we have formulated a mathematical model for the influence of a system of agents. This model qualitatively exhibits behaviours of a system that is guided by an external influencing agent, whereby too weak an influence fails to shape the system and too disparate the external agent from the system causes it to decouple from and ignore the external node. For synchronisation in the Kuramoto model, we are able to detect and quantify these behaviours specifically. We have generated numerical results for the model and shown that an analytical approximation successfully characterises the key behaviours of initial synchronisation between the influencer and the internal system, where the former may be seen as {\it part} of the system and pull their mean-frequency towards itself, to where the influencer is external and stretches or decreases the cohesion of the system, and finally to where the internal system separates and undertakes its own dynamics. Specifically, these results illustrate a parameter regime in which the maximum stretching of the system occurs. The analytical approach also qualitatively reflects these behaviours. Our numerical results demonstrate that even in the deterministic case the time-averaged closeness deviates from zero, illustrating these behaviours through the influence dynamics. The results of this paper show that the region of capture shows a higher-order parameter for the stochastic system, thus demonstrating that the stochastic influence approach is superior to the deterministic one as a representation of social influence. 

Preliminary investigations have provided a foundation for further adaptation of this model. Work in \cite{Gao2014} considers the ability to obtain efficient control via a subset of nodes and to identify the minimum set of so-called driver nodes. We would next like to consider heterogeneous sub-networks for the influence links $B_{i}$. In addition, varying the density of influence nodes such that there is an external network as opposed to a single node is a natural extension for this model. 
To conclude, this paper has demonstrated that mathematical formulations of social influence are possible using non-linear stochastic synchronisation models, distinct from control or the internal states of the latter; where a human agent, `influencer', and the `influencee' are explicitly identified, and the mechanisms of interaction within and across them can be articulated and solved to reveal intuitively natural regimes of behaviour.
\begin{acknowledgments}
This research forms part of a Work-Based PhD at Defence Science and Technology Group through the Australian National University.
\end{acknowledgments}
\nocite{*}

\appendix
\section{\label{app:alphadependence} Dependence on Gamma shape parameter $\alpha$}
\begin{figure}
    \centering
    \includegraphics[width=\linewidth]{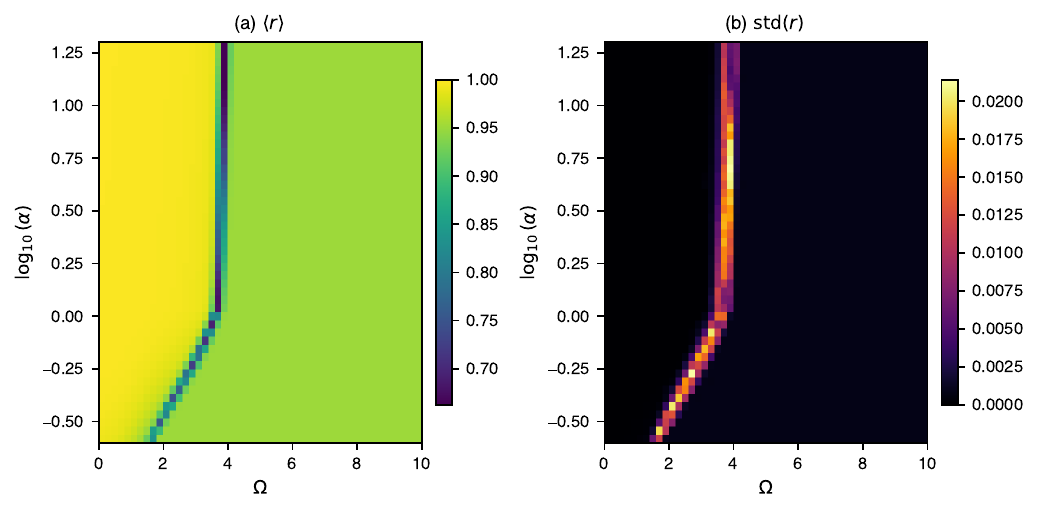}
    \caption{(a) $\langle r \rangle$ in $(\Omega, \log_{10}\alpha)$ space on ER($N{=}40$, $p{=}0.3$) graphs at $\sigma=2$, $\tau = 10$ (100 graph realisations per point, 50 $\Omega$-values in $[0, 10]$, 49 unique $\alpha$-values: 20 log-spaced in $[0.25, 20]$ plus 30 linearly spaced in $[10, 20]$). (b) Corresponding ensemble fluctuations $\mathrm{std}(r)$.}
    \label{fig:alpha_orderparameter}
\end{figure}

Here we consider how the Gamma shape parameter $\alpha$ controls the captured--decoupled landscape, as illustrated in Fig.~\ref{fig:alpha_orderparameter}. We see that for $\alpha>1$ that $\langle r \rangle$ profile is relatively uniform, and thus for this paper we consider only low $\alpha$ such that the Gamma distribution is peaked at $0$. Increasing $\alpha$ shifts the stretching through to larger source frequencies: the minimum moves from $\Omega \approx 1.63$ at $\alpha = 0.25$ to $\Omega\approx 2.45$ at $\alpha = 1$, and to $\Omega \approx 3.7–3.9$ for $\alpha = 10–20$. Thus lower-variance influence links extend the captured regime before the source-frame stretching sets in.As discussed in the introduction, this is representative of an interaction between the external node and the system that is representative of influence. 
 
\section{\label{app:cons}Phase dynamics and conservation law}
\begin{figure}
    \centering
    \includegraphics[width=0.9\columnwidth]{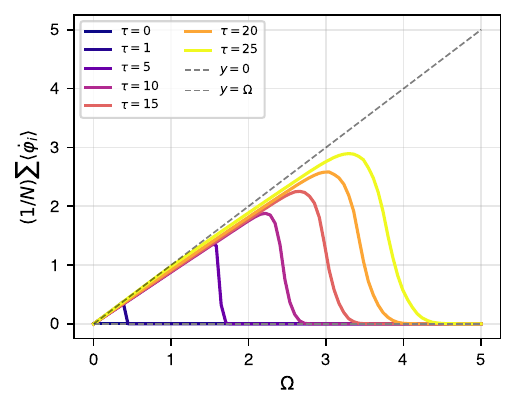}
    \caption{Mean internal oscillator frequency $\left(\frac{1}{N}\right) \sum_{i}\langle\dot{\varphi}_{i}\rangle$ vs $\Omega$ at $\tau \in {0, 1, 5, 10, 15, 20, 25}$ on ER($N=300$, $p=0.3$) graphs ($\sigma = 2$, $\alpha = 1$, 100 graph realisations, 80 $\Omega$-values in [0, 5]). }
    \label{fig:conservation}
\end{figure}

\begin{figure*}
    \centering
    \includegraphics[width=0.97\textwidth]{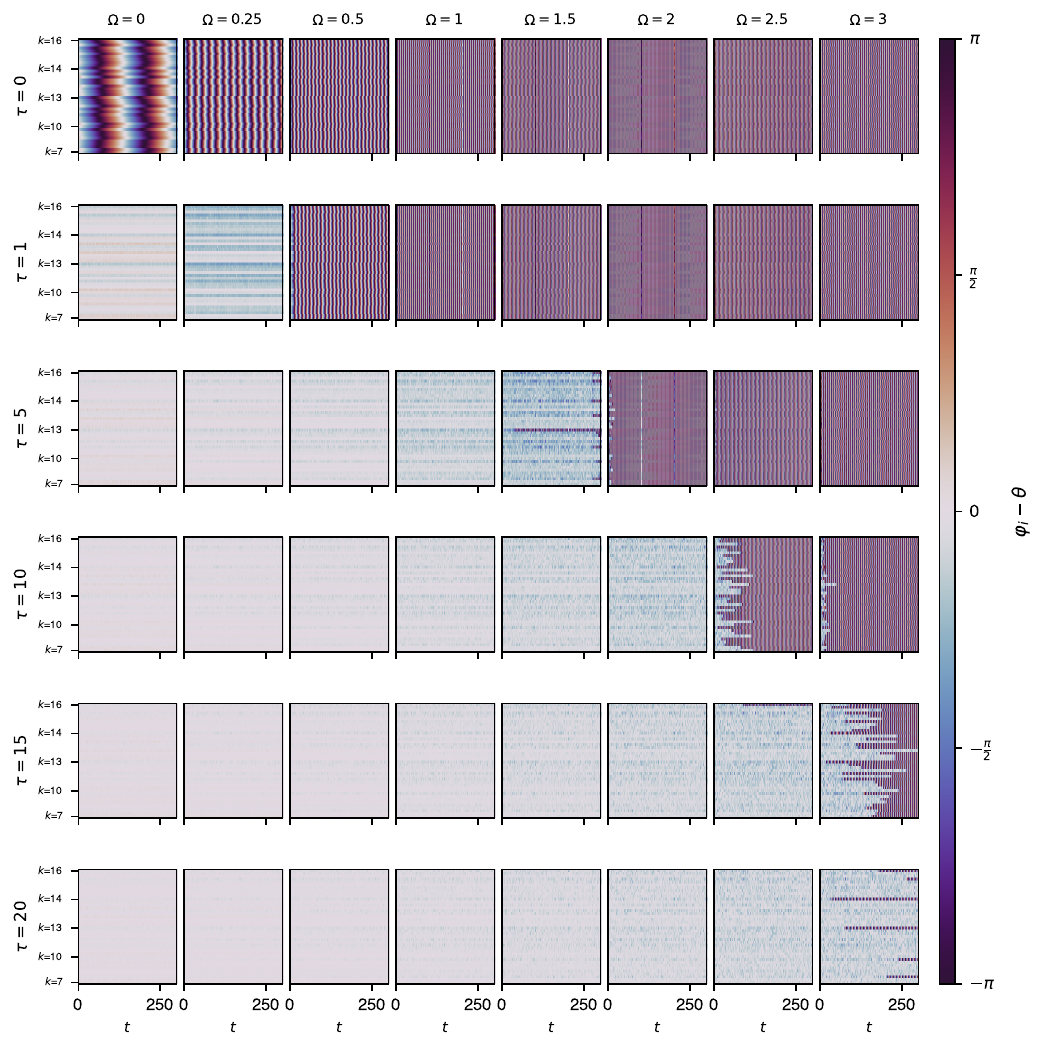}
    \caption{Relative phase $\varphi_i(t) - \theta(t)$ for a single ER($N{=}40$, $p{=}0.3$) graph at $\tau \in \{0, 1, 5, 10, 15, 20\}$, with $\sigma = 2$ and $\alpha=1$. Columns show $\Omega \in \{0, 0.25, 0, 0.25, 0.5, 1.0, 1.5, 2.0, 2.5, 3.0\}$. }
    \label{fig:phase}
\end{figure*}

Fig.~\ref{fig:conservation} looks explicitly at whether the influence term violates the Kuramoto conservation law. In the standard Kuramoto model ($\tau = 0$), the mean frequency remains zero for all $\Omega$, reflecting the conservation law $\frac{d}{dt} \left(\frac{1}{N}\right) \sum_{i} \varphi_{i} = \left(\frac{1}{N}\right) \sum_{i} \omega_{i} = 0$ for the symmetric natural-frequency distribution. With influence, the external node gives the oscillator population a non-zero centre-of-mass frequency over the captured and stretched regimes, before the mean returns to zero in the decoupled high-$\Omega$ regime. The peak mean frequency grows and shifts right with influence strength: approximately 0.37 at $\Omega \approx 0.38$ for $\tau = 1$, 1.37 at $\Omega \approx 1.52$ for $\tau = 5$, 1.88 at $\Omega\approx 2.22$ for $\tau = 10$, 2.58 at $\Omega \approx 3.04$ for $\tau = 20$, and 2.89 at $\Omega \approx 3.29$ for $\tau = 25$. This centre-of-mass diagnostic complements the mean-frequency deviation in Fig.~\ref{fig:regimes} \textbf{Right}: it shows where the external node changes the conserved mean motion, whereas $\Delta$ shows whether individual oscillators are frequency-entrained to the external node.
Fig.~\ref{fig:phase} gives an illustrative example as to how the influence strength shapes the phase dynamics. At $\tau = 0$ the source is disconnected: the oscillator network remains autonomously synchronised, while the source-relative frame winds uniformly as $\Omega$ increases. At $\tau = 1$ only the lowest source frequencies are visibly entrained. Larger $\tau$ progressively extends the captured and stretched regimes: by $\tau$ = 10 coherent source-frame bands persist through $\Omega \approx 2$ and for $\tau = 15–20$ phase slipping is pushed towards the largest displayed source frequencies. The comparison shows that increasing influence strength delays source-frame stretching rather than changing the underlying autonomous synchronisation of the oscillator pack. Comparing rows confirms that $\tau = 0$ satisfies the conservation law $\langle\dot{\varphi}_i\rangle = \omega_i$ (constant relative phases), while increasing $\tau$ progressively violates it.

\section{\label{app:solution}General solutions: fracturing ansatz}
To qualitatively understand the stable solution for this system, we write Eq.~\ref{eq:matrix} as $\dot{X}=\xi-\mathcal{L}X$ and diagonalise $\mathcal{L}$ to
\begin{equation}
\mathcal{L}_D=
  \begin{bmatrix}
   \lambda_{-}   &
   0  \\
   0  &
    \lambda{+}  
   \end{bmatrix}
\end{equation}
such that we now have the form $\dot{X}_{D}=\xi_{D}-\mathcal{L}_{D}X_{D}$, producing two linear differential equations. Using integrating factors, the general solutions are
\begin{align}
X_{D_{-}}&= e^{-\int\lambda_{-}dt}\int\left(e^{\int\lambda_{-}dt}\xi_{D_{-}}dt\right) +C_{-}e^{-\int\lambda_{-}dt}\\
X_{D_{+}}&=e^{-\int\lambda_{+}dt}\int\left(e^{\int\lambda_{+}dt}\xi_{D_{+}}dt\right)+C_{+}e^{-\int\lambda_{+}dt}
\label{eq:Xsol}
\end{align}
where $C_{\pm}$ are integration constants, set by the initial conditions. We note importantly here that there is noise-dependence in both the eigenvalues $\lambda_{\pm}$ through $\bar{G}_1$ and $\bar{G}_2$. We get the following expressions:
\begin{widetext}
\begin{align}
\xi_{D_{-}}&=- \frac{(\Omega-\bar{\omega}_1) \left(\sqrt{4\sigma^2+\tau(\bar{G}_1-\bar{G}_2)^2}-\tau(\bar{G}_1-\bar{G}_2)\right)   -2\sigma(\Omega-\bar{\omega}_2)}{2\sqrt{\sigma^2+2\sigma(\sigma+\tau(\bar{G}_1-\bar{G}_2))+(\sigma-\tau(\bar{G}_1-\bar{G}_2))^2}}\\
\xi_{D_{+}}&=- \frac{(\Omega-\bar{\omega}_1) \left(\sqrt{4\sigma^2+\tau(\bar{G}_1-\bar{G}_2)^2}+\tau(\bar{G}_1-\bar{G}_2)\right)   +2\sigma(\Omega-\bar{\omega}_2)}{2\sqrt{\sigma^2+2\sigma(\sigma+\tau(\bar{G}_1-\bar{G}_2))+(\sigma-\tau(\bar{G}_1-\bar{G}_2))^2}}  \end{align}
\end{widetext}
A qualitative analysis of the solutions $X_{D_{\pm}}$ provides some insights as to the behaviour of the model. We note firstly that the second term $C_{\mp}e^{-\int\lambda_{\mp}dt}\rightarrow 0$ as $t\rightarrow \infty$, thus the solution is dominated by the first term as a time-integral over the random variables $\bar{G}_{1,2}$, which have a time-profile. We also note that for the linearisation above to be valid,  $X_{D}$ must remain \textit{small}. Considering the first term $e^{-\int\lambda_{-}dt}\int\left(e^{\int\lambda_{-}dt}\xi_{D_{-}}dt\right)$ we see that there is a finite window in time, $T_{captured}$, for which $X_{D}$ does remain small. The dependence on $\Omega$ in $\xi_{D_{\mp}}$ suggests that this window decreases as $\Omega$ increases in order for the integral to remain \textit{small}. This affirms the numerical results (see Section \ref{sec:numerical}) which demonstrate that the nodes in the system can be caught to the external node only for a finite a period of time before being dropped. Numerical results also demonstrate this as the system is caught for a decreasing period of time as $\Omega$ increases, as indicated for example in the phase plots shown in Fig.~\ref{fig:phase}. 

When the eigenvalues are time-independent the integrals in Eqs.(\ref{eq:Xsol}) are easily performed with the result
\begin{equation}
    X_{D_\pm} = {\xi_{D_\pm} \over \lambda_{D_\pm}}
    \left( 1- e^{-\lambda_{D_\pm}t}\right) + C_\pm e^{-\lambda_{D_\pm}t} \rightarrow {\xi_{D_\pm} \over \lambda_{D_\pm}}
\label{eq:tindX}
\end{equation}
in the absence of any zero eigenvalue.

Returning to the general expressions above, we note that for $\tau=0$ (where there is no connection to the external node at all) we recover $\xi_{D_{-}}=\frac{1}{2}(\bar{\omega}_1+\bar{\omega}_2)-\Omega$ and $\xi_{D_{+}}=\frac{1}{2}(\bar{\omega}_1-\bar{\omega}_2)$. These correspond to the eigenvalues $\lambda = \{0,2\sigma\} $. This gives the expected result for the stand-alone internal system when fractured, relative to the external node. For the zero mode, however, the solution must be separated out given the singularity in Eq.(\ref{eq:tindX}):
$\dot{X}_{D_-} = \xi_{D_-}$, thus 
$X_{D_-} = \xi_{D_-} t + C_- = \frac{1}{2}\left( (\bar{\omega}_1+\bar{\omega}_2)-\Omega \right)t$ for zero initial conditions.
Thus, the zero mode corresponds to the centre of mass motion expressed in relation to the driver $\Omega$, while the non-zero are stable fluctuations relative to this which will decay for large $t$ leaving only a constant part $X_{D_+} \rightarrow \xi_{D_+}/2\sigma
= \frac{1}{4\sigma}(\bar{\omega}_1-\bar{\omega}_2)$
\cite{Holder2017}. Note also that for equal noise $\bar{G}_1=\bar{G}_2=\bar{G}$
the same solutions for $\xi_{D_{\pm}} $ arise as for $\tau=0$.
Here, both the eigenvalues are non-zero, 
$\lambda=\{\tau \bar{G}, 2\sigma + \tau \bar{G}\}$, and the solutions
$X$ will involve decay.

In the deterministic case, $\bar{G}_1=\bar{G}_2=1$ so that
$\lambda=\{\tau , 2\sigma + \tau \}$ and from Eq.(\ref{eq:tindX}), for large time, we have
\begin{equation}
X_{D} \rightarrow \{ \frac{1}{\tau} \left( \frac{1}{2}(\bar{\omega}_1+\bar{\omega}_2)-\Omega \right) , 
\frac{(\bar{\omega}_1-\bar{\omega}_2)}{2(\tau + 2\sigma)}\}
\end{equation}


\end{document}